\documentclass{article}[12pt]
\usepackage{amsfonts,amsmath,amssymb,latexsym,tikz}
\usepackage{graphicx}
\usepackage{cite}
\newtheorem{thm}{Theorem}[section]
\newtheorem{prp}[thm]{Proposition}

\newtheorem{conj}[thm]{Conjecture}
\newtheorem{cor}[thm]{Corollary}

\newtheorem{claim}{Claim}

\newtheorem{prob}[thm]{Problem}

\newcommand{\qed}{\hfill ~$\square$\bigskip}
\newcommand{\proof}{\noindent{\bf Proof.} }

\newcommand{\cp}{\,\square\,}
\newcommand{\vertex}{\node[vertex]}
\tikzstyle{vertex}=[circle, draw, inner sep=0pt, minimum size=6pt]

\renewcommand{\gg}{\gamma_{g}}
\newcommand{\ggz}{\gamma_{Zg}}
\newcommand{\ggt}{\gamma_{tg}}
\newcommand{\ggl}{\gamma_{Lg}}
\newcommand{\ggll}{\gamma_{LLg}}
\newcommand{\ggs}{\gamma_{g}'}
\newcommand{\ggzs}{\gamma_{Zg}'}
\newcommand{\ggts}{\gamma_{tg}'}
\newcommand{\ggls}{\gamma_{Lg}'}
\newcommand{\gglls}{\gamma_{LLg}'}

\begin{document}

\title{The variety of domination games}

\author{
Bo\v stjan Bre\v sar $^{a,b}$ \and Csilla Bujt{\'a}s $^{c}$ \and
Tanja Gologranc $^{a,b}$
\vspace*{2mm} 
\and Sandi Klav\v zar $^{d,a,b}$ \and Ga\v sper
Ko\v smrlj $^{b,f}$ \and Tilen Marc $^{b,d}$ \vspace*{2mm} 
\and Bal{\'a}zs Patk{\'o}s $^{e}$ \and Zsolt Tuza $^{c,e}$ \and M{\'a}t{\'e} Vizer $^{e}$}

\maketitle

\begin{center}
$^a$ Faculty of Natural Sciences and Mathematics, University of Maribor, Slovenia\\
\medskip

$^b$ Institute of Mathematics, Physics and Mechanics, Ljubljana, Slovenia\\
\medskip

$^c$ Faculty of Information Technology, University of Pannonia, Veszpr\'em, Hungary\\
\medskip

$^d$ Faculty of Mathematics and Physics, University of Ljubljana, Slovenia\\
\medskip

$^e$ Alfr\'ed R\'enyi Institute of Mathematics, Hungarian Academy of
Sciences, Budapest, Hungary
\medskip

$^f$ Abelium R\&D, Ljubljana, Slovenia\\
\medskip

\small \texttt{bostjan.bresar@um.si, \{bujtas,tuza\}@dcs.uni-pannon.hu, sandi.klavzar@fmf.uni-lj.si, tilen.marc@imfm.si, patkos@renyi.hu, gasperk@abelium.eu, \{tanja.gologranc,vizermate\}@gmail.com}
\medskip

\end{center}
\maketitle

\begin{abstract}
Domination game [SIAM J.\ Discrete Math.\ 24 (2010) 979--991] and total domination game [Graphs Combin.\ 31 (2015) 1453--1462] are by now well established games played on graphs by two players, named Dominator and Staller. In this paper, Z-domination game, L-domination game, and LL-domination game are introduced as natural companions of the standard domination games. 

Versions of the Continuation Principle are proved for the new games. It is proved that in each of these games the outcome of the game, which is a corresponding graph invariant, differs by at most one depending whether Dominator or Staller starts the game. The hierarchy of the five domination games is established. The invariants are also bounded with respect to the (total) domination number and to the order of a graph. Values of the three new invariants are determined for paths up to a small constant independent from the length of a path. Several open problems and a conjecture are listed. The latter asserts that the L-domination game number is not greater than $6/7$ of the order of a graph. 
\end{abstract}


\medskip\noindent
{\bf AMS Subj.\ Class.\ (2010)}: 05C69; 05C57; 91A43

\section{Introduction}
\label{sec:intro}

Domination game~\cite{BKR-2010} and total domination game~\cite{HKR-2015} have been investigated in depth by now; see the recent papers~\cite{BDK-2017, BDK-2016, Bu-2015, Bu-2015b, DKR-2015, james-2017, KKS-2016, NSS-2016, XL-2019, XLK-2018} on the domination game,~\cite{HK-2018, HKR-2017, HKR-2018, HR-2017} on the total domination game, as well as references therein. 

In~\cite{bgmrr-2014} the Grundy domination number of a graph $G$ was introduced as the length of a longest sequence of vertices such that each vertex of the sequence dominates at least one new vertex. From our point of view, a vertex by vertex determination of such a sequence can be considered as the domination game for a single player (Staller). The Grundy total domination number was later studied in~\cite{bhr-2016} which can again be considered as a one player total dominaton game. Moreover, motivated by zero forcing sets, Z-Grundy domination number and L-Grundy domination number was investigated in~\cite{BBG-2017}. In this paper we in turn introduce the corresponding Z-domination game, L-domination game, and for reasons to be clarified later, LL-domination game. 

Each of the games, the domination game, the total domination game, and the Z-, L-, and LL-domination game, is played on an isolate-free graph $G$. (Actually, the domination game and the Z-domination game do not require the graph to be isolate-free, but here we will consider only this more restricted case. Nevertheless, to be on the safe side we will state this fact in statements of the results.) All these games can be uniformly described as follows. As usual, for a vertex $v$ of a graph $G$, its open and closed neighborhoods are denoted by $N(v)$ and $N[v]$, respectively. Two players, traditionally named Dominator and Staller, alternately select a vertex from $G$. If Dominator is the first player to select a vertex in a domination game, we speak of a D-game. Otherwise (that is, if Staller begins the game), we have an S-game. In the $i^{\rm th}$ move, the choice of a vertex $v_i$ is legal if for the vertices $v_1,\ldots,v_{i-1}$ chosen so far, the following hold:
    \begin{itemize}
        \item[(i)] $N[v_i] \setminus \bigcup_{j=1}^{i-1}N[v_j]\not=\emptyset$, in
        the domination game;
        \item[(ii)] $N(v_i) \setminus \bigcup_{j=1}^{i-1}N(v_j)\not=\emptyset$, in
        the total domination game;
        \item[(iii)] $N(v_i) \setminus \bigcup_{j=1}^{i-1}N[v_j]\not=\emptyset$, in
        the Z-domination game;
        \item[(iv)] $N[v_i] \setminus \bigcup_{j=1}^{i-1}N(v_j)\not=\emptyset$ and $v_i\neq v_j$ for all $j <i$, in the L-domination game; and
        \item[(v)] $N[v_i] \setminus \bigcup_{j=1}^{i-1}N(v_j)\not=\emptyset$, in
        the LL-domination game.
    \end{itemize}
Each of the games ends if there are no more legal moves available. Dominator wishes to finish the game as soon as possible, while Staller wishes to delay the end. If a D-game is played and both players play optimally, the length of the game, i.e.\ the total number of moves played during the game, is, respectively,  
\begin{itemize}
\item[(i)] the game domination number $\gg(G)$; 
\item[(ii)] the game total domination number $\ggt(G)$; 
\item[(iii)] the game Z-domination number $\ggz(G)$;
\item[(iv)] the game L-domination number  $\ggl(G)$; and
\item[(v)] the game LL-domination number $\ggll(G)$ of $G$. 
\end{itemize}
For the S-game the lengths of the above games give analogous graph invariants $\ggs(G)$, $\ggts(G)$, $\ggzs(G)$, $\ggls(G)$, and $\gglls(G)$, respectively.
 
We proceed as follows. In the next section we prove that the Continuation Principle holds also for the Z-, L-, and LL-domination game. In addition we show that for any of these games the corresponding values of the invariants for the D-game and the S-game differ by at most one. While the proofs of these results for the Z- and LL-domination game are standard (in particular, the difference by at most one follows easily from the corresponding Continuation Principle), the proofs for the L-domination game are more subtle because of the extra condition that the vertices played in the game must be pairwise different. In Section~\ref{sec:hierarchy} we investigate the hierarchy of the five domination games and prove that  $\gg(G)$ and $\ggt(G)$ are upper bounds for $\ggz(G)$ and lower bounds for $\ggl(G)$, which is on the other hand a lower bound for $\ggll(G)$. In Section~\ref{sec:bounds} we bound $\ggz(G)$, $\ggl(G)$ and $\ggll(G)$ and show that  $\ggll(G) \leq n(G)+1$, where $n(G)$ stands for the order of the graph $G$. We also characterize graphs with $\ggll(G)=n(G)+1.$ Then, in Section~\ref{sec:paths}, we consider the path $P_n$ on $n$ vertices and determine the values of $\ggz(P_n)$, $\ggl(P_n)$ and $\ggll(P_n)$ up to an additive constant error term. We conclude the paper with several open problems. In particular, we pose a conjecture that $\ggl(G) \le \frac{6}{7}n(G)$ holds for any isolate-free graph $G$.

\section{Continuation Principles and applications}
\label{sec:continuation}

One of the key tools for the domination game is the Continuation Principle proved in \cite{KWZ-2013}. The corresponding result for the total domination game was established in~\cite{HKR-2015}. Here we prove that the analogous statements, which express the monotonicity of the invariants, are true for the other three domination games. 

For the formulation of the theorem, consider a graph $G$ and a subset $A$ of vertices which are considered to be pre-dominated or pre-totally-dominated. With $G|A$ we denote such a pre-dominated graph, meaning that when a game is played on $G|A$, the vertices from $A$ need not be dominated but they are allowed to be played (provided they are legal moves). More formally, in a Z- or LL-domination game on $G|A$  the choice of a vertex $v_i$ is legal if for the vertices $v_1,\ldots,v_{i-1}$ chosen so far $N(v_i) \setminus (A\cup \bigcup_{j=1}^{i-1}N[v_j])\not=\emptyset$ or $N[v_i] \setminus (A \cup \bigcup_{j=1}^{i-1}N(v_j))\not=\emptyset$ holds, respectively. The condition for the L-domination game on $G|A$ is the same as for the LL-domination game with the additional requirement that no vertex can be selected twice. We use $\ggz(G|A)$, $\ggl(G|A)$, and $\ggll(G|A)$ to denote the number of moves in the Z-, L-, and LL-domination game respectively, under optimal play in a D-game  on $G|A$. For an S-game the analogous invariants are denoted $\ggzs(G|A)$, $\ggls(G|A)$, and $\gglls(G|A)$, respectively.

In the following proofs we will use a standard tool called the imagination strategy that was introduced in the context of the domination game in~\cite{BKR-2010}.

\begin{thm}\label{thm:ContPr}
{\rm (Continuation Principle)}
If $G$ is a graph without isolated vertices and $B\subseteq A \subseteq V(G)$, then
\begin{itemize}
\item[$(i)$] $\ggz(G|A) \le \ggz(G|B)$ and $\ggzs(G|A) \le \ggzs(G|B)$; 
\item[$(ii)$] $\ggll(G|A) \le \ggll(G|B)$ and $\gglls(G|A) \le \gglls(G|B)$; and
\item[$(iii)$] $\ggl(G|A) \le \ggl(G|B)$ and $\ggls(G|A) \le \ggls(G|B)$.
\end{itemize}
\end{thm}

\proof In each proof we consider two games: Game A is the real game played on $G|A$, while Game B is a game on $G|B$ imagined by Dominator. Staller plays optimally in Game A, and Dominator playes optimally in Game B. We denote by $a_i$ and $b_i$ the vertex played in the $i^{\rm th}$ move of Game A and Game B respectively.  If Staller plays $a_i$ in Game A, Dominator copies it into Game B that is $b_i=a_i$ (we will prove that it is always legal). Then, Dominator responds optimally in Game B by playing $b_{i+1}$. If $b_{i+1}$ is a legal move in Game A, Dominator plays the same vertex there i.e, $a_{i+1}=b_{i+1}$. In the other case, $a_{i+1}$ will be defined to be a legal move in Game A (if there exists such a move).  We will prove that the length $\ell_A$ of Game A is not greater than the length $\ell_B$ of Game B.

$(i)$ Consider a real and an imagined Z-domination game as described above. We prove that 
\begin{equation} \label{EQ:1}
B\cup \bigcup_{j=1}^{k}N[b_j] \subseteq A\cup \bigcup_{j=1}^{k}N[a_j]
\end{equation}
 holds for every $k$. Since $B \subseteq A$ is assumed, the analogous relation is valid before the first move ($k=0$). For the inductive step, suppose that  (\ref{EQ:1}) is true for $k=i-1$. If Staller plays $a_i$ in Game A, it is a legal move and hence $N(a_i) \setminus (A\cup \bigcup_{j=1}^{i-1}N[a_j])\not=\emptyset$. Since (\ref{EQ:1}) holds for $k=i-1$, we obtain $N(a_i) \setminus (B\cup \bigcup_{j=1}^{i-1}N[b_j])\not=\emptyset$ that is, $a_i$ is a legal move in Game B. We define $b_i=a_i$ and infer that (\ref{EQ:1}) remains valid with $k=i$. In the other case, the $i^{\rm th}$ move is taken by Dominator. He picks a vertex $b_i$ in Game B. If $b_i$ is a legal move in Game A, he sets $a_i=b_i$ and (\ref{EQ:1}) remains valid for $k=i$. Otherwise, if $b_i$ is not legal in Game A, we have $N(b_i) \subseteq (A\cup \bigcup_{j=1}^{i-1}N[a_j])$ and distinguish two subcases. If $b_i \notin (A\cup \bigcup_{j=1}^{i-1}N[a_j])$, then $a_i$ can be any vertex from  $N(b_i)$, and (\ref{EQ:1}) remains valid. If $b_i \in (A\cup \bigcup_{j=1}^{i-1}N[a_j])$, we have $B\cup \bigcup_{j=1}^{i}N[b_j] \subseteq A\cup \bigcup_{j=1}^{i-1}N[a_j]$. Consequently, $a_i$ can be any legal move in Game A, (\ref{EQ:1}) remains valid for $k=i$. 

This proves that (\ref{EQ:1}) holds after each move. Then, Game A cannot be longer then Game B that is $\ell_A \le \ell_B$. Since Staller played optimally on $G|A$ and Dominator played optimally on $G|B$, we may conclude that $\ggz(G|A)\le \ell_A \le \ell_B \le \ggz(G|B)$. Our proof equally holds for the D-game and the S-game. Thus $\ggzs(G|A)\le \ell_A \le \ell_B \le \ggzs(G|B)$ also follows.

$(ii)$ Now, consider a real and an imagined LL-domination game on $G|A$ and $G|B$ respectively. Here, we prove  
\begin{equation} \label{EQ:2}
B\cup \bigcup_{j=1}^{k}N(b_j) \subseteq A\cup \bigcup_{j=1}^{k}N(a_j)
\end{equation}
 for every $k$. 
The argumentation is very similar to that of part $(i)$. The assumption $B\subseteq A$ gives the base for the inductive proof. If Staller playes $a_i$ in Game A and (\ref{EQ:2}) holds with $k=i-1$, then $N[a_i] \setminus (B\cup \bigcup_{j=1}^{i-1}N(b_j))\supseteq N[a_i] \setminus (A\cup \bigcup_{j=1}^{i-1}N(a_j)) \not= \emptyset$. Therefore $b_i=a_i$ is a legal move in Game B and, after this choice, (\ref{EQ:2}) is valid for $k=i$. Now, suppose that Dominator plays the $i^{\rm th}$ move $b_i$ in Game B. If $b_i$ is a legal move in Game A, we set $a_i=b_i$ and (\ref{EQ:2}) clearly holds with $k=i$. If $b_i$ is not a legal move in the real game,  $N[b_i] \subseteq A\cup \bigcup_{j=1}^{i-1}N(a_j)$. Then, already the set $A\cup \bigcup_{j=1}^{i-1}N(a_j)$ is a superset of $B\cup \bigcup_{j=1}^{i}N(b_j)$ and for any legal move $a_i$, the relation (\ref{EQ:2}) will be valid for $k=i$. Since Game A finishes when $A\cup \bigcup_{j=1}^{k}N(a_j) =V(G)$, the conclusion is $\ell_A\le \ell_B$ and the required inequalities follow as in the proof of $(i)$.

$(iii)$ In an L-domination game, no vertex can be played more than once. Hence, we define the following sets $F_k^A$ and $F_k^B$ containing those vertices which cannot be played after the $k^{\rm th}$ move in the real and in the imagined L-domination game, respectively:
\begin{eqnarray}
\label{F1}
F_k^A  =  \{a_1, \dots , a_k\} \cup \{v \in V(G): N[v] \subseteq A \cup \bigcup_{j=1}^{k}N(a_j)\},\\
\label{F2}
F_k^B  =  \{b_1, \dots , b_k\} \cup \{v \in V(G): N[v] \subseteq B \cup \bigcup_{j=1}^{k}N(b_j)\}. 
\end{eqnarray}
Observe that for any two sets $S$, $S'$ of vertices, $S' \subseteq S $ implies $\{v \in V(G): N[v] \subseteq S'\} \subseteq \{v \in V(G): N[v] \subseteq S\}$. We will prove that both (\ref{EQ:2}) and
\begin{equation} \label{EQ:3}
F_k^B \subseteq F_k^A
\end{equation} 
hold for every $k \ge 0$. Our condition $B \subseteq A$ implies that (\ref{EQ:2}) and (\ref{EQ:3}) are valid for $k=0$. For the inductive proof, assume that both (\ref{EQ:2}) and (\ref{EQ:3}) hold for $k=i-1$. If Staller plays $a_i$ in Game A, then $N[a_i] \setminus (A\cup \bigcup_{j=1}^{i-1}N(a_j))\neq \emptyset$ and, by (\ref{EQ:2}), $N[a_i] \setminus (B\cup \bigcup_{j=1}^{i-1}N(b_j))\neq \emptyset$ also holds.  Moreover, by (\ref{EQ:3}), $a_i\notin F_{i-1}^A$ implies $a_i \notin F_{i-1}^B$. Then, $b_i=a_i$ is a legal move in Game B. This definition ensures that (\ref{EQ:2}) and (\ref{EQ:3}) hold for $k=i$.
In the other case Dominator plays $b_i$ as the $i^{\rm th}$ move in Game B. If $b_i$ is also a legal move in Game A, we set $a_i=b_i$. Then both~\eqref{EQ:2} and~\eqref{EQ:3} hold with $k=i$. Now, assume that $b_i$ is not a legal move in the real game that is $b_i \in F_{i-1}^A$. Then, at least one of the following statements is true:
\begin{itemize}
\item[$(a)$] $b_i=a_s$ for some $s<i$;
\item[$(b)$] $N[b_i] \subseteq A \cup \bigcup_{j=1}^{i-1}N(a_j)$.
\end{itemize}
In either case, we have $N(b_i) \subseteq (A \cup \bigcup_{j=1}^{i-1}N(a_j))$ and hence, 
$$B\cup \bigcup_{j=1}^{i}N(b_j) \subseteq A\cup \bigcup_{j=1}^{i-1}N(a_j).$$ The latter relation and $b_i \in F_{i-1}^A$ also imply $F_i^B \subseteq F_{i-1}^A$. Then, $a_i$ can be chosen as any legal move in Game A, (\ref{EQ:2}) and (\ref{EQ:3}) will be valid with $k=i$. 

We have just proved that (\ref{EQ:2}) can be maintained for the real and imagined games. This implies, as in the previous parts of the proof, that $\ell_A \leq \ell_B$ and the two inequalities stated in $(iii)$ follow.  
\qed

Next we prove that for the Z-, L- and LL-domination games, the lengths of the D-game and S-game cannot differ by more than 1, if the players play optimally. The analogous statements for domination and total domination game were proved in \cite{BKR-2010, HKR-2015, KWZ-2013}.

\begin{thm} 
\label{thm:SD}
If $G$ is a graph without isolated vertices, then 
\begin{itemize}
\item[$(i)$] $|\ggz(G)-\ggzs(G)| \le 1$; 
\item[$(ii)$] $|\ggll(G)-\gglls(G)| \le 1$; and 
\item[$(iii)$] $|\ggl(G)-\ggls(G)| \le 1$. 
\end{itemize} 
\end{thm}
\proof
For the Z- and LL-domination games, the statements easily follow from Theorem~\ref{thm:ContPr} $(i)$ and $(ii)$. First, consider a Z-domination game on $G$ where Staller starts by playing one of her optimal first moves, say $a_1$. If the game is not finished yet, from the second move we may interpret it as a D-game on $G|N[a_1]$ and therefore, by Theorem~\ref{thm:ContPr} $(i)$, we have
$$ \ggzs(G) = 1+ \ggz(G|N[a_1]) \le 1+\ggz(G).$$
Similarly, if Dominator starts the game on $G$ and $b_1$ is one of his optimal first moves,
$$ \ggz(G) = 1+ \ggzs(G|N[b_1]) \le 1+\ggzs(G)$$
follows. These establish part $(i)$.

The proof is very similar for the LL-domination game, so we omit it. The main difference is that here $G|N(a_1)$ and $G|N(b_1)$ are considered instead of $G|N[a_1]$ and $G|N[b_1]$.

To prove $(iii)$, we cannot use Theorem~\ref{thm:ContPr} $(iii)$ directly, but we can use the imagination strategy. First, assume that Staller starts the L-domination game on $G$  by playing $a_1$ which is an optimal first move. This will be the real game, Game A, where Staller plays optimally. The imagined game (i.e., Game B) will be an L-domination game on $G$ where Dominator starts by playing $b_1$ and he plays optimally throughout. Hence, the move $a_1$ is not copied into Game B. On the other hand, for every odd $i$ with $i\ge 3$, the move $a_i$ will be copied into Game B by setting  $b_{i-1}=a_i$. Dominator chooses his moves to be optimal in Game B, and for every positive odd $i$, his move $b_i$ is copied into Game A as $a_{i+1}=b_i$, if it is legal in the real game; otherwise we show that $a_{i+1}$ can be any legal move in Game A. Our main goal is to prove that 
\begin{equation} \label{EQ:4}
\bigcup_{j=1}^{k-1}N(b_j) \subseteq \bigcup_{j=1}^{k}N(a_j)
\end{equation}
holds for all $k\ge 1$. We also define the set of forbidden vertices $F_k^A$ and $F_k^B$ after the $k^{\rm th}$ move, similarly as they were given in~\eqref{F1} and~\eqref{F2}, respectively, that is, 
\begin{eqnarray*}
F_k^A  =  \{a_1, \dots , a_k\} \cup \{v \in V(G): N[v] \subseteq \bigcup_{j=1}^{k}N(a_j)\},\\
F_k^B  =  \{b_1, \dots , b_k\} \cup \{v \in V(G): N[v] \subseteq \bigcup_{j=1}^{k}N(b_j)\}. 
\end{eqnarray*}
We are going to prove that
\begin{equation} \label{EQ:5}
F_{k-1}^B \subseteq F_k^A
\end{equation} 
holds for every $k\ge 1$.

Relations (\ref{EQ:4}) and (\ref{EQ:5}) clearly hold for $k=1$, since in this case the left-hand side sets are empty sets. For the inductive step, we suppose and (\ref{EQ:4}) and (\ref{EQ:5}) are true for $k=i-1\ge 0$. If $i$ is odd, $a_i$ is chosen optimally by Staller in Game A. Since it is a legal move, $a_i \notin F_{i-1}^A$, and the inductive hypothesis on (\ref{EQ:5}) implies $a_i \notin F_{i-2}^B$. Hence, $b_{i-1}=a_i$ is a legal move in Game B and we may conclude that both (\ref{EQ:4}) and (\ref{EQ:5}) remain valid with $k=i$. If $i$ is even, $b_{i-1}$ is chosen optimally by Dominator in Game B. If it is legal in Game A, we set $a_i=b_{i-1}$ and then,  (\ref{EQ:4}) and (\ref{EQ:5}) hold with $k=i$. If $b_{i-1}$ is not legal in Game A that is, $b_{i-1} \in F_{i-1}^A$, we have the following possibilities:
\begin{itemize}
\item[$(a)$] $b_{i-1}=a_s$ for some $s<i$; or
\item[$(b)$] $N[b_{i-1}] \subseteq \bigcup_{j=1}^{i-1}N(a_j)$.
\end{itemize}
In either case, we have $N(b_{i-1}) \subseteq  \bigcup_{j=1}^{i-1}N(a_j)$ that, together with  (\ref{EQ:4}), implies
$$\bigcup_{j=1}^{i-1}N(b_j) \subseteq \bigcup_{j=1}^{i-1}N(a_j).$$ 
By this relation, by the inductive hypothesis on (\ref{EQ:5}), and since $b_{i-1} \in F_{i-1}^A$, we conclude than $F_{i-1}^B \subseteq F_{i-1}^A$.
Consequently, choosing an arbitrary legal move  $a_{i}$ in Game A, (\ref{EQ:4}) and (\ref{EQ:5}) remain valid with $k=i$. Therefore, we can maintain (\ref{EQ:4}) and (\ref{EQ:5}) during the games while letting  Staller and Dominator play optimally in Game A and B respectively. For the lengths $\ell_A$ and $\ell_B$ of Game A and B, (\ref{EQ:4}) implies $\ell_A-1 \le \ell_B$. Finally, we get 
$$\ggls(G) \le \ell_A \le \ell_B+1 \le \ggl(G)+1.$$

To prove $\ggl(G) \le \ggls(G)+1$, we can use an analogous argumentation. Here the real game, namely Game A, is an L-domination game where Dominator starts with an optimal first move $a_1$ and then, from the second turn, Staller plays optimally. Game B is the imagined L-domination game in which Dominator plays optimally. Moreover, this is an S-game which begins with copying here the move $a_2$ from Game A as $b_1$. The moves are defined by the rules used in the previous proof. More precisely, if $i$ is even, Staller picks an optimal move  $a_i$ in Game A, and it can be proved that $b_{i-1}=a_i$ is a legal move in Game B. If $i$ is an odd number with $i \ge 3$, Dominator plays an optimal vertex $b_{i-1}$ in Game  B and we define $a_i=b_{i-1}$ in Game A, if it is  legal; otherwise, we may define $a_i$ to be an arbitrary but legal move in Game A.
As it can be proved, again,  (\ref{EQ:4}) and (\ref{EQ:5}) hold for every $k\ge 1$. These imply $\ell_A  \le \ell_B+1$. Since Staller plays optimally in Game A and Dominator in Game B, we have
$$\ggl(G) \le \ell_A \le \ell_B+1 \le \ggls(G)+1$$
which completes the proof of the theorem.
\qed

\section{Hierarchy of the games}
\label{sec:hierarchy}

In this section we show that the five domination games fulfill the following hierarchy. 

\begin{thm}
\label{thm:hierarchy}
If $G$ is a graph without isolated vertices, then
$$\ggz(G)\le \gg(G), \ggt(G) \le \ggl(G)\le \ggll(G)\,.$$
\end{thm}

\proof To prove $\ggz(G)\le \gg(G)$, consider a domination game on $G$ and assume
that Dominator plays optimally. Suppose further that in the $i$th
move he selects a vertex $v_i$ such that $N[v_i] \setminus
\bigcup_{j=1}^{i-1}N[v_j]=\{v_i\}$. If he chooses a neighbor $u \in
N(v_i)$ instead of $v_i$,
$\bigcup_{j=1}^{i-1}N[v_j] \cup N[u] \supseteq
\bigcup_{j=1}^{i}N[v_j]$ holds, and by the Continuation Principle for the domination game (
\cite{KWZ-2013}) this change does not lengthen the game. Hence,
Dominator can choose a strategy in the domination game which is optimal and also obeys the rule
$N(v_i) \setminus \bigcup_{j=1}^{i-1}N[v_j]\neq \emptyset$ of the
Z-domination game in each of his turns. Therefore, it means real restriction only for Staller when Z-domination game is played instead of
domination game on $G$. Thus, $\ggz(G)\le \gg(G)$.

The remaining inequalities will be proved by using the imagination strategy.
Given an isolate-free graph $G$, we assume that Dominator and Staller play a real
game on $G$, while Staller also imagines another game is played on it.
Dominator plays optimally in the real game and Staller plays
optimally in the imagined one. We denote by $v_i$ and
$u_i$ the vertex played in the $i^{\rm th}$ move of the real and the
imagined game respectively, and prove that Staller can guarantee  that the length $\ell'$ of the imagined game
is not greater than the length $\ell$ of the real game.

\medskip
We next prove that $\ggz(G)\le \ggt(G)$. Here the real game is a total domination
game on $G$ and the imagined one played by Staller is a  Z-domination game on $G$.
We will prove that
\begin{equation}
\label{eq:Game1}
\bigcup_{j=1}^{i}N[u_j] \supseteq \bigcup_{j=1}^{i}N(v_j)
\end{equation}
holds    for every $i$. Since the Z-domination game ends when
$\bigcup_{j=1}^{i}N[u_j]=V(G)$ and the total domination game ends when
$\bigcup_{j=1}^{i}N(v_j)=V(G)$, the condition~\eqref{eq:Game1} will imply $\ell' \le \ell$.

In the first turn, Dominator plays $v_1$ in the real game and
Staller copies it into the imagined game ($u_1=v_1$). Clearly,
(\ref{eq:Game1}) is valid for $i=1$.
Now, assume that (\ref{eq:Game1}) is true for $i=k-1$. If $k$ is
even, Staller selects a vertex $u_k$ in the Z-domination game.
This move is also legal in the real game, because $N(u_k) \setminus \bigcup_{j=1}^{k-1}N[u_j]\neq \emptyset$
and (\ref{eq:Game1}) together imply that $N(u_k) \setminus
\bigcup_{j=1}^{k-1}N(v_j)\neq \emptyset$. Then, Staller copies
her move into the real game ($v_k=u_k$) and (\ref{eq:Game1}) will be valid for $i=k$. In the other case
$k$ is odd and Dominator plays a vertex $v_k$ in the total domination game.
If $u_k=v_k$ is a legal move in the imagined game, Staller just copies
the move into the Z-domination game and (\ref{eq:Game1}) remains valid for
$i=k$. If $v_k$ is not a legal choice in the Z-domination game,
$N(v_k) \subseteq \bigcup_{j=1}^{k-1}N[u_j]$ and hence, after any
legal choice of $u_k$, the relation (\ref{eq:Game1}) remains
valid.
Hence,  $\ell' \le \ell$. Moreover,
since Staller played optimally in the imagined Z-domination game and Dominator
played optimally in the total domination game, we have $\ggz(G) \le \ell'
\le \ell \le \ggt$.

\medskip
We proceed with the proof of $ \gg(G) \le \ggl(G)$. Now, the real game is an L-domination game and the imagined one is
a domination game on $G$. We prove that~\eqref{eq:Game1}
holds for all $i$. In the first turn, Dominator's move $v_1$ is
copied into the imagined game and (\ref{eq:Game1}) clearly holds
for $i=1$. Suppose that (\ref{eq:Game1}) is valid for $i=k-1 \ge
1$. Under this assumption, if $k$ is even, every legal (and optimal) move  $u_k$ of Staller in the
imagined game is legal in the real game as well. Setting
$v_k=u_k$,~\eqref{eq:Game1} will be valid with $i=k$.  If $k$ is
odd, any optimal move $v_k$ of Dominator in the real game is
either valid in the imagined game and $u_k=v_k$ maintains the
relation (\ref{eq:Game1}), or
$N[v_k] \subseteq \bigcup_{j=1}^{k-1}N[u_j]$ holds. In the latter
case, from our assumption $\bigcup_{j=1}^{k-1}N[u_j] \supseteq \bigcup_{j=1}^{k}N(v_j)$ follows, that
any legal move $u_k$ in the imagined game ensures that
(\ref{eq:Game1}) is valid for $i=k$. Thus, the imagined game
finishes no later than the real one and we have
$\gg(G) \le \ell' \le \ell \le \ggl(G)$.

\medskip
Next we prove $\ggt(G) \le \ggl(G)$. To prove this inequality, the real game is an L-domination game and
the imagined one is a total domination game on $G$. As we will see,
Staller may ensure that
\begin{equation}
\label{eq:Game2}
\bigcup_{j=1}^{i}N(u_j) \supseteq \bigcup_{j=1}^{i}N(v_j)
\end{equation}
holds  for every $i$.
Setting $u_1=v_1$, (\ref{eq:Game2}) will be valid for $i=1$.
Then, assume that (\ref{eq:Game2}) is valid for $i=k-1$. If $k$
is even, any legal move $u_k$ taken in the total domination game is
also legal in the L-domination game, and (\ref{eq:Game2}) remains valid for
$i=k$, if $u_k$ is copied into the real game. The situation is similar
if $k$ is odd and $v_k$ is a legal move in the imagined game.
The only remaining case is when Dominator chooses a vertex $v_k$ in the
L-domination game such that $N(v_k)\subseteq \bigcup_{j=1}^{k-1}N(u_j)$.
Then, again, any legal move $u_k$ in the total domination game ensures that
(\ref{eq:Game2}) is valid for $i=k$. These prove
$\ggt(G) \le \ell' \le \ell \le \ggl(G)$.

\medskip
It remains to prove that $\ggl(G)\le \ggll(G)$. Let the real game be the LL-domination game and let the imagined one
be the L-domination game. Staller plays by applying an optimal strategy in the L-domination game and
meantime, she ensures that (\ref{eq:Game2}) is valid after each
turn. Assuming that (\ref{eq:Game2}) holds with $i=k-1$, any move $u_k$ of
Staller may be copied into the real game and (\ref{eq:Game2})
remains valid. If $k$ is odd and the move $v_k$ of Dominator is
not legal in the L-domination game, we have either $ \bigcup_{j=1}^{k-1}N(u_j)\supseteq N[v_k]$ or $v_k=u_j$ for a $j<k$.
In both cases, our assumption implies
$ \bigcup_{j=1}^{k-1}N(u_j)\supseteq \bigcup_{j=1}^{k}N(v_j)$.
Therefore, any legal choice of $u_k$ maintains relation
(\ref{eq:Game2}). At the end, we obtain $\ggl(G) \le \ell' \le \ell \le \ggll(G)$.
\qed

\begin{figure}[ht!]
    \begin{center}
        \begin{tikzpicture}[]
        \tikzstyle{vertex}=[circle, draw, inner sep=0pt, minimum size=2pt]
        \tikzset{vertexStyle/.append style={rectangle}}
        \vertex (1) at (0,0) [label=below:$\ggz$] {};
        \vertex (2) at (-1,1) [label=left:$\ggt$] {};
        \vertex (3) at (1,1) [label=right:$\gg$] {};
        \vertex (4) at (0,2) [label=left:$\ggl \;$] {};
        \vertex (5) at (0,2.8) [label=left:$\ggll$] {};
        \path
        (1) edge[dashed] (2)
        (1) edge[dashed] (3)
        (2) edge[dashed] (4)
        (3) edge[dashed] (4)
        (4) edge[dashed] (5);
        \end{tikzpicture}
    \end{center}
    \caption{Relations between the five versions of the game domination
        number.} 
\label{fig:Hasse}
\end{figure}
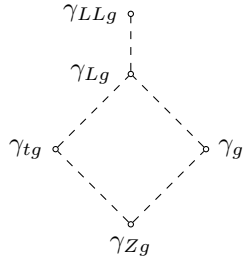

Theorem~\ref{thm:hierarchy} together with fact (established in~\cite{HKR-2015}) that $\gg$ and $\ggt$ are incomparable, can be briefly presented with the Hasse diagram representing the partial
ordering  between $\ggz$, $\gg$, $\ggt$,
$\ggl$, and $\ggll$ as shown in Fig.~\ref{fig:Hasse}. 

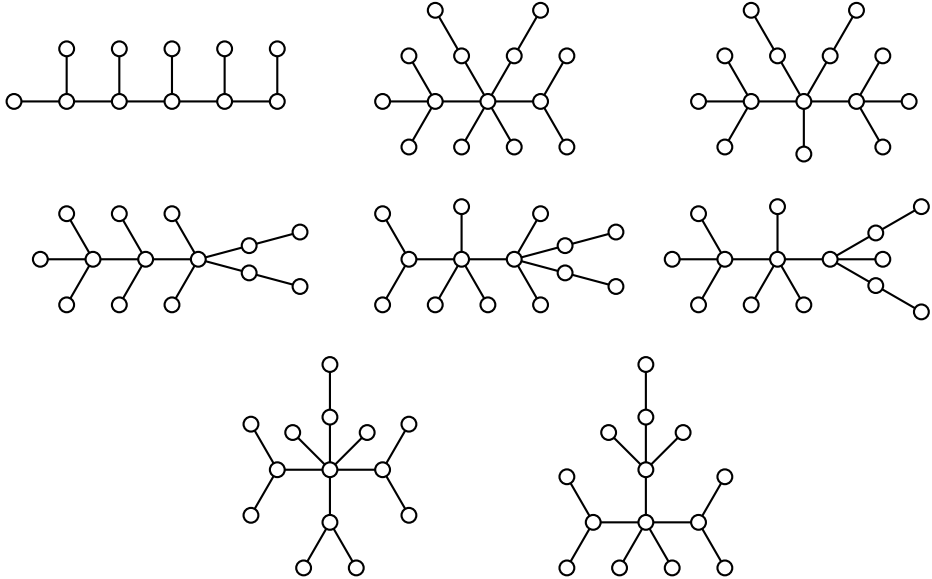
\begin{figure}[ht!]
\begin{center}
\begin{tikzpicture}[scale=0.7,style=thick]
\def\vr{2pt/0.5}

\coordinate(g11x1y1) at (-0.5,-3);
\begin{scope}[shift=(g11x1y1)]
\foreach \i in {0,...,5}
{
	\coordinate(x\i) at (\i,0);
}

\coordinate(x6) at (5,1);

\draw(x0) -- (x1) -- (x2) -- (x3) -- (x4) -- (x5) -- (x6);

\coordinate(x7) at ({1+1*cos(90.0)}, {1*sin(90.0)});
\draw(x1) -- (x7);

\coordinate(x8) at ({2+1*cos(90.0)}, {1*sin(90.0)});
\draw(x2) -- (x8);

\coordinate(x9) at ({3+1*cos(90.0)}, {1*sin(90.0)});
\draw(x3) -- (x9);

\coordinate(x10) at ({4+1*cos(90.0)}, {1*sin(90.0)});
\draw(x4) -- (x10);

\foreach \i in {0,...,10}
{
	\draw (x\i)[fill=white] circle (\vr);
}
\end{scope}

\coordinate(g14x1y2) at (6.5,-3);
\begin{scope}[shift=(g14x1y2)]
\foreach \i in {0,...,3}
{
	\coordinate(x\i) at (\i,0);
}

\coordinate(x4) at ({3+1*cos(60.0)}, {1*sin(60.0)});

\draw(x0) -- (x1) -- (x2) -- (x3) -- (x4);

\coordinate(x5) at ({1+1*cos(-120.0)}, {1*sin(-120.0)});
\draw(x1) -- (x5);

\coordinate(x6) at ({1+1*cos(120.0)}, {1*sin(120.0)});
\draw(x1) -- (x6);

\coordinate(x7) at ({2+1*cos(-120.0)}, {1*sin(-120.0)});
\draw(x2) -- (x7);

\coordinate(x8) at ({2+1*cos(60.0)}, {1*sin(60.0)});
\coordinate(x9) at ({2+2*cos(60.0)}, {2*sin(60.0)});
\draw(x2) -- (x8) -- (x9);

\coordinate(x10) at ({2+1*cos(120.0)}, {1*sin(120.0)});
\coordinate(x11) at ({2+2*cos(120.0)}, {2*sin(120.0)});
\draw(x2) -- (x10) -- (x11);

\coordinate(x12) at ({2+1*cos(-60.0)}, {1*sin(-60.0)});
\draw(x2) -- (x12);

\coordinate(x13) at ({3+1*cos(-60.0)}, {1*sin(-60.0)});
\draw(x3) -- (x13);

\foreach \i in {0,...,13}
{
	\draw (x\i)[fill=white] circle (\vr);
}
\end{scope}

\coordinate(g14x1y3) at (12.5,-3);
\begin{scope}[shift=(g14x1y3)]
\foreach \i in {0,...,4}
{
	\coordinate(x\i) at (\i,0);
}

\draw(x0) -- (x1) -- (x2) -- (x3) -- (x4);

\coordinate(x5) at ({1+1*cos(-120.0)}, {1*sin(-120.0)});
\draw(x1) -- (x5);

\coordinate(x6) at ({1+1*cos(120.0)}, {1*sin(120.0)});
\draw(x1) -- (x6);

\coordinate(x7) at ({2+1*cos(60.0)}, {1*sin(60.0)});
\coordinate(x8) at ({2+2*cos(60.0)}, {2*sin(60.0)});
\draw(x2) -- (x7) -- (x8);

\coordinate(x9) at ({2+1*cos(-90.0)}, {1*sin(-90.0)});
\draw(x2) -- (x9);

\coordinate(x10) at ({2+1*cos(120.0)}, {1*sin(120.0)});
\coordinate(x11) at ({2+2*cos(120.0)}, {2*sin(120.0)});
\draw(x2) -- (x10) -- (x11);

\coordinate(x12) at ({3+1*cos(60.0)}, {1*sin(60.0)});
\draw(x3) -- (x12);

\coordinate(x13) at ({3+1*cos(-60.0)}, {1*sin(-60.0)});
\draw(x3) -- (x13);

\foreach \i in {0,...,13}
{
	\draw (x\i)[fill=white] circle (\vr);
}
\end{scope}

\coordinate(g14x2y1) at (0.0,-6);
\begin{scope}[shift=(g14x2y1)]
\foreach \i in {0,...,3}
{
	\coordinate(x\i) at (\i,0);
}

\draw(x0) -- (x1) -- (x2) -- (x3);

\coordinate(x6) at ({1+1*cos(-120.0)}, {1*sin(-120.0)});
\draw(x1) -- (x6);

\coordinate(x7) at ({1+1*cos(120.0)}, {1*sin(120.0)});
\draw(x1) -- (x7);

\coordinate(x8) at ({2+1*cos(-120.0)}, {1*sin(-120.0)});
\draw(x2) -- (x8);

\coordinate(x9) at ({2+1*cos(120.0)}, {1*sin(120.0)});
\draw(x2) -- (x9);

\coordinate(x10) at ({3+1*cos(120.0)}, {1*sin(120.0)});
\draw(x3) -- (x10);

\coordinate(x11) at ({3+1*cos(15.0)}, {1*sin(15.0)});
\coordinate(x12) at ({3+2*cos(15.0)}, {2*sin(15.0)});
\draw(x3) -- (x11) -- (x12);

\coordinate(x4) at ({3+1*cos(-15.0)}, {1*sin(-15.0)});
\coordinate(x5) at ({3+2*cos(-15.0)}, {2*sin(-15.0)});
\draw(x3) -- (x4) -- (x5);

\coordinate(x13) at ({3+1*cos(-120.0)}, {1*sin(-120.0)});
\draw(x3) -- (x13);

\foreach \i in {0,...,13}
{
	\draw (x\i)[fill=white] circle (\vr);
}
\end{scope}

\coordinate(g14x2y2) at (6.0,-6);
\begin{scope}[shift=(g14x2y2)]
\foreach \i in {1,...,3}
{
	\coordinate(x\i) at (\i,0);
}

\draw (x1) -- (x2) -- (x3);

\coordinate(x0) at ({1+1*cos(-120.0)}, {1*sin(-120.0)});
\draw(x1) -- (x0);

\coordinate(x6) at ({1+1*cos(120.0)}, {1*sin(120.0)});
\draw(x1) -- (x6);

\coordinate(x7) at ({2+1*cos(-60.0)}, {1*sin(-60.0)});
\draw(x2) -- (x7);

\coordinate(x8) at ({2+1*cos(90.0)}, {1*sin(90.0)});
\draw(x2) -- (x8);

\coordinate(x9) at ({2+1*cos(-120.0)}, {1*sin(-120.0)});
\draw(x2) -- (x9);

\coordinate(x10) at ({3+1*cos(60.0)}, {1*sin(60.0)});
\draw(x3) -- (x10);

\coordinate(x11) at ({3+1*cos(15.0)}, {1*sin(15.0)});
\coordinate(x12) at ({3+2*cos(15.0)}, {2*sin(15.0)});
\draw(x3) -- (x11) -- (x12);

\coordinate(x4) at ({3+1*cos(-15.0)}, {1*sin(-15.0)});
\coordinate(x5) at ({3+2*cos(-15.0)}, {2*sin(-15.0)});
\draw(x3) -- (x4) -- (x5);

\coordinate(x13) at ({3+1*cos(-60.0)}, {1*sin(-60.0)});
\draw(x3) -- (x13);

\foreach \i in {0,...,13}
{
	\draw (x\i)[fill=white] circle (\vr);
}
\end{scope}

\coordinate(g14x2y3) at (12.0,-6);
\begin{scope}[shift=(g14x2y3)]
\foreach \i in {0,...,4}
{
	\coordinate(x\i) at (\i,0);
}

\draw(x0) -- (x1) -- (x2) -- (x3) -- (x4);

\coordinate(x6) at ({1+1*cos(-120.0)}, {1*sin(-120.0)});
\draw(x1) -- (x6);

\coordinate(x7) at ({1+1*cos(120.0)}, {1*sin(120.0)});
\draw(x1) -- (x7);

\coordinate(x8) at ({2+1*cos(-60.0)}, {1*sin(-60.0)});
\draw(x2) -- (x8);

\coordinate(x9) at ({2+1*cos(90.0)}, {1*sin(90.0)});
\draw(x2) -- (x9);

\coordinate(x10) at ({2+1*cos(-120.0)}, {1*sin(-120.0)});
\draw(x2) -- (x10);

\coordinate(x11) at ({3+1*cos(30.0)}, {1*sin(30.0)});
\coordinate(x5) at ({3+2*cos(30.0)}, {2*sin(30.0)});
\draw(x3) -- (x11) -- (x5);

\coordinate(x12) at ({3+1*cos(-30.0)}, {1*sin(-30.0)});
\coordinate(x13) at ({3+2*cos(-30.0)}, {2*sin(-30.0)});
\draw(x3) -- (x12) -- (x13);

\foreach \i in {0,...,13}
{
	\draw (x\i)[fill=white] circle (\vr);
}
\end{scope}

\coordinate(g14x3y1) at (4.5,-10);
\begin{scope}[shift=(g14x3y1)]
\foreach \i in {0,...,2}
{
	\coordinate(x\i) at (\i,0);
}

\draw(x0) -- (x1) -- (x2);

\coordinate(x3) at ({1*cos(120.0)}, {1*sin(120.0)});
\draw(x0)--(x3);

\coordinate(x4) at ({1*cos(-120.0)}, {1*sin(-120.0)});
\draw(x0)--(x4);

\coordinate(x5) at ({2+1*cos(60.0)}, {1*sin(60.0)});
\draw(x2)--(x5);

\coordinate(x6) at ({2+1*cos(-60.0)}, {1*sin(-60.0)});
\draw(x2)--(x6);

\coordinate(x7) at (1, -1);
\draw(x1)--(x7);

\coordinate(x8) at ({1+cos(-60.0)}, {-1+sin(-60.0)});
\draw(x7) -- (x8);

\coordinate(x9) at ({1+cos(-120.0)}, {-1+sin(-120.0)});
\draw(x7) -- (x9);

\coordinate(x10) at ({1+1*cos(45.0)}, {1*sin(45.0)});
\draw(x1) -- (x10);

\coordinate(x11) at ({1+1*cos(135.0)}, {1*sin(135.0)});
\draw(x1) -- (x11);

\coordinate(x12) at (1, 1);
\coordinate(x13) at (1, 2);
\draw(x1) -- (x12) -- (x13);

\foreach \i in {0,...,13}
{
	\draw (x\i)[fill=white] circle (\vr);
}
\end{scope}

\coordinate(g14x3y2) at (10.5,-11);
\begin{scope}[shift=(g14x3y2)]
\foreach \i in {0,...,2}
{
	\coordinate(x\i) at (\i,0);
}

\draw(x0) -- (x1) -- (x2);

\coordinate(x3) at ({1*cos(120.0)}, {1*sin(120.0)});
\draw(x0)--(x3);

\coordinate(x4) at ({1*cos(-120.0)}, {1*sin(-120.0)});
\draw(x0)--(x4);

\coordinate(x5) at ({2+1*cos(60.0)}, {1*sin(60.0)});
\draw(x2)--(x5);

\coordinate(x6) at ({2+1*cos(-60.0)}, {1*sin(-60.0)});
\draw(x2)--(x6);

\coordinate(x7) at (1, 1);
\draw(x1)--(x7);

\coordinate(x8) at ({1+cos(45.0)}, {1+sin(45.0)});
\draw(x7) -- (x8);

\coordinate(x9) at ({1+cos(135.0)}, {1+sin(135.0)});
\draw(x7) -- (x9);

\coordinate(x12) at (1, 2);
\coordinate(x13) at (1, 3);
\draw(x7) -- (x12) -- (x13);

\coordinate(x10) at ({1+1*cos(-60.0)}, {1*sin(-60.0)});
\draw(x1) -- (x10);

\coordinate(x11) at ({1+1*cos(-120.0)}, {1*sin(-120.0)});
\draw(x1) -- (x11);

\foreach \i in {0,...,13}
{
	\draw (x\i)[fill=white] circle (\vr);
}
\end{scope}

\end{tikzpicture}
\end{center}
\caption{Smallest trees with pairwise different values of the game domination numbers.}
\label{fig:smallest-trees-with-different-values}
\end{figure}

The five game domination invariants from Theorem~\ref{thm:hierarchy} can be pairwise different as we have demonstrated with a computer search over the class of trees. The smallest such trees are presented in Fig.~\ref{fig:smallest-trees-with-different-values}.

The top-left tree on $11$ vertices in Fig.~\ref{fig:smallest-trees-with-different-values} has the following values: 
$$\ggz = 5, \gg=6, \ggt=7, \ggl=8, \ggll=9\,.$$ 
The remaining trees are the smallest trees (each of them having $14$ vertices) with the same separability property, except that $\gg$ and $\ggt$ are reversed. More precisely, for these seven trees we have
$$\ggz = 5, \ggt=6, \gg=7, \ggl=8, \ggll=9\,.$$

Note that all five inequalities in
Theorem~\ref{thm:hierarchy} are sharp. For example, for the
five-cycle, $\ggz(C_5)= \gg(C_5)= \ggt(C_5) =
\ggl(C_5) =3$ (but $\ggll(C_5)=5$). Another example,
say $F_n$, can be obtained from the complete graph $K_n$ by
attaching $n$ leaves to each vertex of $K_n$. Assume that $n \ge 2$.
Clearly, $\gamma_t(F_n)=n$. First consider the Z-domination game and the total
domination game on $F_n$ and assume that Dominator plays a vertex
$v_1\in V(K_n)$ in the first turn. Then, in the Z-domination game, we have
$N(u)\setminus N[v_1]=\emptyset$ for every leaf $u$. Hence, all the
played vertices are from $V(K_n)$ and $\ggz(F_n)=n$. In the
total domination game, if Dominator plays $v_1$ as his first move,
Staller may respond by playing a leaf adjacent to $v_1$. But once at
least two non-leaf vertices have been played in the total domination
game, no leaves can be chosen in the later turns. Hence, Staller can
ensure that at least one leaf is played, and Dominator has a
strategy to ensure that at most one leaf is played. This implies
$\ggt(F_n)=n+1$. In the domination game, L- and LL-domination game,
Staller always may play a leaf while there is at least one
non-selected vertex of higher degree. Hence, $\gg(F_n) \ge
2n-1$. On the other hand, $2\gamma_t(F_n)-1=2n-1 \ge
\ggll(F_n)$ (see Proposition~\ref{prp:bounds} below) implies $\gg(F_n)= \ggl(F_n)=
\ggll(F_n)=2n-1$.

\section{Bounds on $\ggz$, $\ggl$, and $\ggll$}
\label{sec:bounds}

In this section we bound $\ggz(G)$, $\ggl(G)$, and $\ggll(G)$. In the main result we prove that $\ggll(G)\le n(G)+1$ and characterize the equality case. In this way we round off Theorem~\ref{thm:hierarchy}. 

To state the results we need to recall some standard terminology. We say that a vertex of a graph {\em totally dominates} its neighbors and  {\em dominates} itself and its neighbors. If $S$ is a subset of vertices of a graph $G$, then $S$ {\em (totally) dominates} $G$ if each vertex of $G$ is (totally) dominated by some vertex from $S$. The size of a smallest set that (totally) dominates $G$ is called {\em (total) domination number} of $G$. These invariants are denoted by $\gamma(G)$ (resp.\ $\gamma_t(G)$).  

The bounds $\gamma(G)\le \gg(G)\le 2\gamma(G)-1$ and $\gamma_t(G)\le \ggt(G)\le 2\gamma_t(G)-1$ were proved in~\cite{BKR-2010} and~\cite{HKR-2015}, respectively. For the three domination games introduced in this paper the following related bounds hold. 

\begin{prp} 
\label{prp:bounds}
If $G$ is a graph without isolated vertices, then
    \begin{itemize}
        \item[$(i)$] $\gamma(G)\le \ggz(G)\le 2\gamma(G)-1$;
        \item[$(ii)$] $\gamma_t(G)\le \ggl(G)\le 2\gamma_t(G)-1$; and
        \item[$(iii)$] $\gamma_t(G)+1\le \ggll(G)\le
        2\gamma_t(G)-1$.
    \end{itemize}
\end{prp}

\proof 
A Z-domination game ends when the set $\{v_1,\ldots, v_i\}$ of the
chosen vertices becomes a dominating set of $G$. Indeed, otherwise
we have a vertex $u$ outside of $\bigcup_{j=1}^{i}N[v_j]$ and the choice of a
neighbor of $u$ is a legal move. Hence, $\gamma(G)\le
\ggz(G)$. The upper bound follows from Theorem~\ref{thm:hierarchy} and the above mentioned inequality $\gg(G)\le 2\gamma(G)-1$. This proves (i). 

The L-domination game and the
LL-domination game on $G$ ends when $\bigcup_{j=1}^{i}N(v_j)=V(G)$, moreover
Dominator may fix a total dominating set $D'$ and in each move he
plays a vertex $w\in D'$  for which
$N(w)\setminus \bigcup_{j=1}^{i}N(v_j)\neq \emptyset$ (while such a vertex exists). This implies the upper bounds in (ii) and (iii). The lower bound in (ii) is a direct consequence of Theorem~\ref{thm:hierarchy} and the inequality $\gamma_t(G) \le \ggt(G)$. 
Concerning the lower bound in $(iii)$, we note that $2 \le \gamma_t (G) \le \ggll(G)$ holds for every graph $G$. Moreover, in the second move of an  LL-domination game Staller may repeat the first move of Dominator. In this way she can ensure that $\ggll(G) \ge \gamma_t(G) + 1$.
\qed

In the rest of this section we prove an upper bound on $\ggll(G)$, for which we need the following two results. 

\begin{prp}\label{prp:llparity}
If $G$ is a graph without isolated vertices and Staller plays an LL-domination game according to an optimal strategy, then Dominator makes the last move of the game. In particular, $\ggll(G)$ is odd and $\gglls(G)$ is even. In particular, $\ggll(G)\neq \gglls(G)$.
\end{prp}

\proof
Suppose $u\in V(G)$ is the last move in an LL-domination game on $G$ with Staller playing optimally. Then there exists $v\in N[u]\setminus \cup_{i=1}^m N(u_i)$ with $u_1, \ldots, u_m$ being the vertices picked earlier in the game. We infer that $v\ne u$ because otherwise $u$ could be selected again, so the game would not be over yet. Hence Staller was not the player who has selected $u$, because she could play $v$ instead of $u$ and prolong the game for at least one more move. 
\qed
 
One can strengthen Proposition \ref{prp:llparity} using the Continuation Principle.

\begin{prp}\label{prp:llcomponent}
If $G$ is a graph without isolated vertices and an LL-domination game is played on $G$, then there exists an optimal strategy of Staller such that the last move in \textit{every component} is made by Dominator.
\end{prp}

\proof
Suppose Staller, according to an optimal strategy $\mathcal{S}$, plays the last move $u$ of a component $C$ during an LL-domination game on a graph $G$ with previous moves $u_1,\ldots, u_m$. Then there exists $v\in N[u]\setminus \cup_{i=1}^mN(u_i)$. If $v=u$, then $u$ is still a legal move, so Staller does not finish the game on $C$. If $v\neq u$, then Staller could have played $v$ instead of $u$. After her move, $v$ would still be a legal move, while if $u$ is played, no vertex in $C$ is legal. So by Theorem~\ref{thm:ContPr}(ii), the latter strategy is at least as good as $\mathcal{S}$.
\qed

All is now ready for the main result of this section. 

\begin{thm}\label{thm:llbound}
If $G$ is a graph without isolated vertices, then $\ggll(G)\le n(G)+1$. Moreover, equality holds if and only if all components of $G$ are $K_2$s.
\end{thm}

\proof
First we prove the theorem for connected graphs. We start with a simple claim.

\begin{claim}\label{claim:llskipvertex}
If the minimum degree in $G$ is at least $2$, then in any LL-domination game on $G$ at least one vertex will not be picked at all.
\end{claim}

\noindent{\bf Proof (of Claim 1).}
Observe that at the moment when the last but first vertex is picked all vertices are totally dominated.
\qed

Suppose the vertices picked by the players are $u_1,\ldots, u_m$. Let $d_m$ denote the number of repetitions, i.e. $|\{j\le m: \exists i<j ~ u_i=u_j\}|$. Furthermore, let $b_m$ denote the number of isolated vertices in $G[u_1,\ldots,u_m]$.

\begin{claim}\label{claim:lldomstrat}
If Dominator starts the game, then he has a strategy that for any $m$ we have $d_{2m+1}+b_{2m+1}\le 1$ and if Staller starts the game, then he can manage to maintain $d_{2m}+b_{2m}=0$.
\end{claim}

\noindent{\bf Proof (of Claim 2).}
We proceed by induction on $m$. If Staller starts by picking $u_1$, then Dominator picks a neighbor $u_2$ of $u_1$ and thus $d_2=0, b_2=0$. If Dominator starts, then he can pick $u_1$ arbitrarily. If Staller picks $u_2=u_1$, then Dominator picks $u_3\in N(u_1)$ and obtains $d_3=1,b_3=0$. If Staller picks $u_2\neq u_1$, then Dominator picks $u_3\in N(u_1)\cup N(u_2)\setminus \{u_1,u_2\}$ (as $G$ is connected, if no such vertex exists, then there is no legal move for Dominator) and obtains $d_3=0,b_3\le 1$.

For the inductive step in the Staller starts game: first of all, as $b_{2m}=0$, Staller cannot pick a previously played vertex. If Staller picks a vertex $u_{2m+1}\notin \cup_{j=1}^{2m}N(u_j)$, then Dominator picks any neighbor $u_{2m+2}$ of $u_{2m+1}$. If $u_{2m+1}\in \cup_{j=1}^{2m}N(u_j)$, then Dominator picks any legal vertex $u_{2m} \in \cup_{j=1}^{2m+1}N(u_j)$. Note that as $G$ is connected, if no legal vertex exists, then the game is over. In both cases, Dominator maintained $d_{2m+2}=b_{2m+2}=0$.

For the inductive step in the Dominator starts game: first of all, Staller can pick a previously played vertex only if $b_{2m+1}=1$ and thus $d_{2m+1}=0$ and Staller can only repeat the isolated vertex of $G[u_1,\ldots,u_{2m+1}]$. So if $u_{2m+2}$ is a repeated vertex, then Dominator can pick any neighbor $u_{2m+3}$ of $u_{2m+2}$ obtaining $d_{2m+3}=1,b_{2m+3}=0$. If $b_{2m+1}=0$, then Dominator proceeds as in the Staller-start-game.
\qed
\smallskip

Clearly, Claim \ref{claim:lldomstrat} proves $\ggll(G)\le n+1$ and $\gglls(G)\le n$ for any connected graph $G$ on $n$ vertices. Also, Claim \ref{claim:llskipvertex} yields $\ggll(G)\le n$ for any connected graph $G$ with minimum degree at least two. Suppose $G$ contains a vertex $v$ of degree one. Then Dominator modifies his strategy as follows: he first picks a neighbor $u_1$ of $v$. Then depending on Staller's move $u_2$, he responds as follows: 
\begin{itemize}
\item
If $u_2=u_1$, then he picks $u_3\in N(u_1)\setminus \{v\}$ (note that $u_3$ exists if $G$ is not $K_2$). This ensures that $v$ cannot be picked during the game. At this point, we have $d_3=1,b_3=0$ and Dominator is able to follow his strategy above to guarantee $\ggll(G)\le n$.
\item
If $u_2\in N(u_1)$, then $d_2=b_2=0$ and the above strategy of Dominator guarantees $\ggll(G)\le n$.
\item
Finally, if $u_2\notin N(u_1)$, then Dominator picks $u_3\in N(u_1)\setminus \{v\}$ to ensure that $v$ cannot be picked during the game. At this point, we have $d_3=0,b_3=1$ and Dominator is able to follow his strategy above to guarantee $\ggll(G)\le n$.
\end{itemize}
This concludes the proof if $G$ is connected.

For the general case let $G$ be an isolate-free graph on $n$ vertices with at least one component $C_1$ consisting of at least 3 edges. We can assume that Staller follows a strategy as in Proposition \ref{prp:llcomponent}, so in each component, Dominator makes the last move. Then Dominator starts by picking a vertex $u_1\in C_1$ according his strategy above for $C_1$. By the assumption on Staller's strategy, Dominator can always play a vertex from the component of Staller's last move and therefore partition the game into games on the components in such a way that all components' games are Staller-start-games apart from the one on $C_1$. So $\ggll(G)\le \ggll(C_1)+\sum_{i=2}^k\gglls(C_i)\le n$, where $C_2,\ldots, C_k$ are the other components of $G$.
\qed

Combining Theorem~\ref{thm:llbound} and Proposition~\ref{prp:bounds}(iii) with Theorem~\ref{thm:hierarchy} we have: 

\begin{cor}
\label{cor:hierarchy}
If $G$ is a graph without isolated vertices, then
$$\ggz(G)\le \gg(G), \ggt(G) \le \ggl(G)\le \ggll(G) \le \min\{2\gamma_t(G)-1, n(G)+1\}\,.$$
\end{cor}

\section{The games played on paths}
\label{sec:paths}

In this section we examine the values of the five games on one of the simplest graphs, the path graphs. The result for the domination game was first proved in the unpublished manuscript~\cite{KWZ2012}, an alternative proof appeared several years later in~\cite{Ko-2017}. It reads as follows: 

\begin{thm}
\label{thm:dg-path}
{\em \cite{KWZ2012, Ko-2017}}  If $n\ge 1$, then
$$\gamma_g(P_n) = 
\left\{
\begin{array}{ll}
\left\lceil\frac{n}{2}\right\rceil-1; & n\equiv 3\pmod 4\,, \\ \\
\left\lceil\frac{n}{2}\right\rceil; & otherwise\,.
\end{array}
\right.$$
\end{thm}

Dorbec and Henning~\cite{dh-2016} obtained the corresponding result for the total domination game. 

\begin{thm}
\label{thm:tdg-path}
{\em \cite{dh-2016}}  If $n\ge 2$, then
$$\ggt(P_n) = 
\left\{
\begin{array}{ll}
\left\lfloor \frac{2n}{3} \right\rfloor; & n\equiv 5\pmod 6\,, \\ \\
\left\lceil\frac{2n}{3}\right\rceil; & otherwise\,.
\end{array}
\right.$$
\end{thm}

\begin{figure}[ht!]
    \begin{center}
        \begin{tikzpicture}[]
        \tikzstyle{vertex}=[circle, draw, inner sep=0pt, minimum size=2pt]
        \tikzset{vertexStyle/.append style={rectangle}}
        \vertex (1) at (0,0) [label=below:$\ggz(P_n)\approx \frac{n}{2}$] {};
        \vertex (2) at (-1,1) [label=left:$\ggt(P_n)\approx \frac{2n}{3}$] {};
        \vertex (3) at (1,1) [label=right:$\gg(P_n)\approx \frac{n}{2}$] {};
        \vertex (4) at (0,2) [label=left:$\ggl(P_n)\approx \frac{2n}{3}\;$] {};
        \vertex (5) at (0,2.8) [label=left:$\ggll(P_n)\approx \frac{4n}{5}$] {};
        \path
        (1) edge[dashed] (2)
        (1) edge[dashed] (3)
        (2) edge[dashed] (4)
        (3) edge[dashed] (4)
        (4) edge[dashed] (5);
        \end{tikzpicture}
    \end{center}
    \caption{The five domination games played on the path graphs.} 
\label{fig:games-on-paths}
\end{figure}
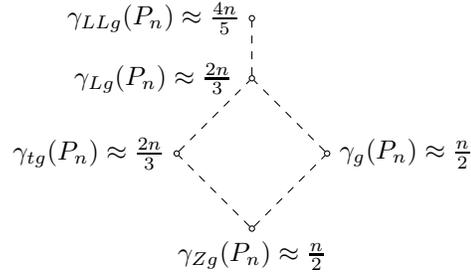

Hence $\gg(P_n)$ roughly equals $n/2$, while $\ggt (P_n)$ is roughly $2n/3$. In the following we will prove similar results for the other three games, where we will consider only approximate values since obtaining the exact ones could double the length of proofs with tedious case analysis. The asymptotics of the parameters for all of the five games is presented in Fig.~\ref{fig:games-on-paths}. 

For the rest of the section we assume that $V(P_n) = \{0, 1,\ldots, n-1\}$, where the vertices appear in the natural order, that is, $i$ and $j$ are connected by an edge if and only if $|i-j|=1$.

\begin{thm}
For every positive integer $n$ there exists a constant $c_n$ such that $\ggz(P_n) = \frac{n}{2} + c_n$ holds with $|c_n|\le 2$.
\end{thm}

\proof
The upper bound  is obtained by Theorems~\ref{thm:hierarchy} and \ref{thm:dg-path}, i.e., $\ggz (P_n)\leq \gg (P_n) \le \frac{n}{2} + \frac{1}{2}$.

To obtain a lower bound we will consider a strategy for Staller. Suppose the $i^{\rm th}$ move is Staller's move and let vertex $k$ be the smallest vertex not in $\bigcup_{j=1}^{i-1}N[v_j]$. Then the vertex $k-1$ is a legal move with $|\bigcup_{j=1}^{i}N[v_j]| - |\bigcup_{j=1}^{i-1}N[v_j]| = 1$, unless $k=0$. In the latter case, vertex 1 is a legal move only needed at most once at Staller's first move. On the other hand, at each move Dominator can dominate at most three new vertices. Hence, besides possibly thefirst move, Staller has a strategy to achieve that in every two moves at most four new vertices are dominated. This implies a lower bound $\ggz (P_n)\geq  2+\frac{n-6}{2}-1=\frac{n}{2}-2$.
\qed

\begin{thm}
For every positive integer $n$ there exists a constant $c_n$ such that $\ggl(P_n) = \frac{2n}{3} + c_n$ holds with $|c_n|\le 1$.
\end{thm}

\proof
The lower bound is obtained by Theorems \ref{thm:hierarchy} and \ref{thm:tdg-path}, i.e.~$\ggl(P_n)\geq \ggt(P_n) \ge \frac{2n}{3} -1$.

To obtain an upper bound we will provide a strategy for Dominator. For every $i$ we define three values:
\begin{itemize}
\item Let $p_i= | \bigcup_{j=1}^{i}N(v_j)|$.
\item Consider graphs $G^1$ and $G^2$ whose vertices are $\{0, 2,\ldots, 2 \lfloor \frac{n}{2}\rfloor\}$ and $\{1, 3,\ldots, 2 \lceil \frac{n}{2}\rceil - 1\}$ and are both isomorphic to paths with the increasing order of vertices. Let $G_i^1$ and $G_i^2$ denote the induced subgraphs of $G^1$ and $G^2$ on the vertices in $\bigcup_{j=1}^{i}N(v_j)$. Then let  $d_i$ be the total number of connected components in $G_i^1$ and $G_i^2$. The empty graph has one connected component.
\item Let $f_i$ denote the number of vertices $v$ after the $i^{\rm th}$ move such that $v$ was not chosen before and is not in $\bigcup_{j=1}^{i}N(v_j)$ but $N(v)$ is a subset of $\bigcup_{j=1}^{i}N(v_j)$.
\end{itemize}

The strategy of Dominator is the following: say that after the $i^{\rm th}$ move vertex $k$ is the smallest vertex not in $\bigcup_{j=1}^{i}N(v_j)$. Then Dominator chooses $v_{i+1}$ to be $k+1$ and thus totally dominates $k$ and maybe also $k+2$. We claim that with such a move it holds $\Delta := (p_{i+1} - d_{i+1} -f_{i+1})-(p_i - d_i -f_i) \geq 2$. Notice that $d_{i+1} - d_i\leq 0$ by the choice of the move, unless $i+1=3$, $v_3=1$ and $v_2$ is an odd ineger picked by Staller. First assume that $p_{i+1} - p_{i} = 2$. Observe that $f_{i+1}-f_{i}$ could be positive only because of the vertices $k-1$ or $k+3$, if they were not chosen before, but after this move their open neighborhoods are totally dominated. But $k$ is the smallest vertex not totally dominated, thus $k-1$ cannot increase $f_i$. But if the open neighborhood of $k+3$  is totally dominated after the $(i+1)^{\rm st}$ move, then $k+4$ was totally dominated after the $i^{\rm th}$ move, implying that $d_{i+1}-d_{i}= -1$. Thus in this case $\Delta = 2$. If $f_{i+1}-f_{i}$ is not positive, then clearly $\Delta \geq 2$. Now assume that $p_{i+1} - p_{i} = 1$. This implies that $k+3$ was chosen before and that $d_{i+1}-d_{i}= -1$. Hence  $k+3$ cannot increase $f_i$, while $k-1$ cannot do it by the same reasoning as before. Thus in this case $\Delta \geq 2$.

In the final part of the proof we show that on the Staller's move $\Delta \geq 1$. First case is if Staller plays a move when it does not increase $p_i$. Then $d_i$ remains the same, while $f_i$ decreases by one giving $\Delta = 1$. If Staller plays a move when $p_{i+1} - p_{i} = 1$, it is easy to see that then either $d_{i+1}-d_{i}=f_{i+1}-f_{i}=0$ or $d_{i+1}-d_{i} = -1$ and $f_{i+1}-f_{i} \leq 1$ giving $\Delta \geq 1$. Finally if $p_{i+1} - p_{i} = 2$, then only one connected component can be created, but in this case $f_{i+1}-f_{i}=0$.  On the other hand, $f_{i+1}-f_{i}=2$ implies $d_{i+1}-d_{i}=-1$, thus also in this case  $\Delta \geq 1$.

We have proved that with this strategy of Dominator  for every $i$ we have $(p_{i+2} - d_{i+2} -f_{i+2})-(p_i - d_i -f_i) \geq 3$ unless $i+2$ is 3 or 4 and even then $(p_{i+2} - d_{i+2} -f_{i+2})-(p_i - d_i -f_i) \geq 2$. Hence $p_{m} - d_{m} -f_{m} +2 \geq \frac{3m}{2}-1$, where $m$ is the final length of the game. But $p_m = n$, $d_m = 2$, and $f_m =0$, thus $m \leq \frac{2n}{3}+1$. 
\qed

In the proof of the following theorem we will consider also predominated graphs. If $G$ is a graph and $v$ a vertex of $G$, then we say that $v$ is predominated for the LL-domination game if the move $v_i$, for which $N[v_i] \setminus \bigcup_{j=1}^{i-1}N(v_j) = \{v\}$, is forbidden.

\begin{thm}
It holds $\ggll(P_n) = \frac{4n}{5} + c_n$ for some small bounded constants $c_n$.
\end{thm}

\proof
As above let $p_i= | \bigcup_{j=1}^{i}N(v_j)|$.
First we present a strategy for Staller showing a lower bound for $\ggll(P_n)$. The strategy is the following: if Staller can play a move for which $p_i$ does not increase, then this move is played. Otherwise, if vertex $k$ is the smallest vertex not totally dominated, then Staller chooses $k-1$ (or $k+1$ if $k=0$ and thus $i=2$). We prove that if $i\ge 4$ and the $i^{\rm th}$ move is played by Staller, the $p_{i+3}-p_{i-1}\le 5$, i.e., the value of $p_i$ increases by at most 5 within four consecutive moves. Notice that by the choice of the Staller's moves, each Staller's move can increase $p_i$ by at most $1$. Also, by the definition of the game, $p_i$ can increase by at most $2$ on Dominator's turn. Hence we must prove that in four moves it cannot happen that Dominator increases $p_i$ twice by $2$ and Staller twice by $1$. 

Assume that Staller played $v_i$ with $p_i-p_{i-1}=1$ and Dominator picked $v_{i+1}=k$ with $p_{i+1}-p_i=2$. Note that this implies $k-v_i\ge 3$ and $v_j\neq k,k-2,k+2$ for all $j<i$. We claim that Staller can at the $(i+2)^{\rm nd}$ move repeat the vertex $k$. Assume that this is not the case, i.e., vertex $k$ is already totally dominated, i.e. for some $j<i$ we have $v_j=k-1$ or $k+1$. But then could have selected $k-1$ or $k+1$ for the $i^{\rm th}$ move without increasing $p_i$, contrary to the assumption. Thus vertex $k$ is not totally dominated at the $(i+2)^{\rm nd}$ move. Thus if $m$ is the total number of moves during the play, then $n\le p_m\le 7+\frac{5}{4}(m-4)$ showing $\ggll (P_n)\geq  \frac{4n}{5} -2$.

To prove the upper bound we will consider a Staller start LL-domination game on two predominated graphs. Let $P_n^1$ be the predominated graph $P_n$ with vertices 0 and $n-1$ predominated. Similarly let $P_n^2$ be the predominated graph $P_n$ with vertices 0, 2, 4, and $n-1$ predominated. We will prove that:
\begin{equation}
\gglls (P_n^1) \leq \theta(P_n^1) :=
\left\{
\begin{array}{ll}
4 \left\lfloor \frac{n}{5} \right\rfloor; & n\equiv 0,1,2 \pmod 5\,, \\ \\
4 \left\lceil\frac{n}{5}\right\rceil + 2; & n\equiv 3,4 \pmod 5\,,
\end{array}
\right.
\end{equation}
and that
\begin{equation}
\gglls(P_n^2) \leq \theta(P_n^2) :=
\left\{
\begin{array}{ll}
4 \left\lfloor \frac{n}{5} \right\rfloor -2; & n\equiv 0,1 \pmod 5\,, \\ \\
4 \left\lfloor \frac{n}{5} \right\rfloor; & n\equiv 2,3 \pmod 5\,, \\ \\
4 \left\lfloor \frac{n}{5} \right\rfloor +2; & n\equiv 4 \pmod 5\,.
\end{array}
\right.
\end{equation}

We prove the above statements by induction on $n$ for $n\geq 3$ in the case of $P_n^1$, and for $n\geq 5$ in the case of $P_n^2$. We have calculated the exact numbers of $\gglls(P_n^1)$ and $\gglls(P_n^2)$ up to $n=24$ and the above values are exact.

Consider $P_n^1$ with $n$ reasonably big (say $n>24$), and assume that the inductive assumption holds for  $P_m^1,P_m^2$, with $m<n$. We consider various first moves of Staller.
\begin{itemize}
\item Staller chooses $v_1=1$ or $v_1=2$: then Dominator can choose $v_2=2$ in the first case and $v_1=1$ in the second. The game on the obtained predominated graph with predominated vertices 0, 1, 2, 3, and $n-1$ is equivalent to the game on $P_{n-3}^1$.  It holds $\gglls(P_n^1) \leq \gglls(P_{n-3}^1) + 2$ which is at most $\theta(P_n^1)$ in all the cases of $n \text{ (mod 5)}$.
\item Staller chooses $v_1=n-2$ or $v_1=n-3$: symmetric as above.
\item Staller chooses $v_1=0$: then Dominator can choose $v_2=4$. The game on the obtained predominated graph with predominated vertices 0, 1, 3, 5, and $n-1$ is equivalent to the game on $P_{n-1}^2$.  It holds $\gglls(P_n^1) \leq \gglls(P_{n-1}^2) + 2$ which is exactly $\theta(P_n^1)$ in all the cases of $n \text{ (mod 5)}$.
\item Staller chooses $v_1=n-1$: symmetric as above.
\item Staller chooses $2<v_1<n-3$: If Dominator chooses either $v_2=v_1-1$ or $v_2=v_1+1$, then the obtained predominated graph has predominated vertices either $0,v_1-2, v_1-1, v_1, v_1+1, n-1$ or $0,v_1-1, v_1, v_1+1, v_1+2, n-1$. In particular, Dominator can consider the same strategy as playing on two disjoint graphs: $P_{v_1}^1$ and $P_{n-(v_1+2)}^1$ in the first case, or $P_{v_1-1}^1$ or $P_{n-(v_1+1)}^1$ in the second. We have:
\begin{eqnarray*}
\gglls(P_n^1) & \leq & \min\{\theta (P_{v_1}^1) + \theta (P_{n-(v_1+2)}^1) + 2, \theta (P_{v_1-1}^1) + \theta (P_{n-(v_1+1)}^1) + 2\} \\
& \leq & \theta (P_n^1)\,.
\end{eqnarray*}
In fact the last inequality holds since $\theta (P_n^1) < \theta (P_{v_1}^1) + \theta (P_{n-(v_1+2)}^1) + 2$ only if $v_1 = 0 \text{ (mod 5)}$ and $n-(v_1+2) = 0 \text{ (mod 5)}$ as it can be checked by an easy examination. In that case the second entry of the minimum is smaller. 
\end{itemize}

Now consider $P_n^2$ with $n$ reasonably big (say $n>24$) and again assume that the inductive assumption holds for  $P_m^1,P_m^2$, with $m<n$. Similarly as above we consider various first moves of Staller:
\begin{itemize}
\item Staller chooses $v_1=1$ or $v_1=3$: then Dominator can select $v_2=2$. The game on the obtained predominated graph with predominated vertices 0, 1, 2, 3, 4 and $n-1$ is equivalent to the game on $P_{n-4}^1$.  It holds $\gglls(P_n^2) \leq \gglls(P_{n-4}^1) + 2$ which is exactly $\theta(P_n^1)$ in all the cases of $n \text{ (mod 5)}$.
\item Staller chooses $v_1=0$, $v_1=2$, or $v_1=4$: then Dominator can choose $v_2=4$ in the first two cases, and $v_1=2$ in the third. The game on the obtained predominated graph with predominated vertices 0, 1, 2, 3, 4, 5 and $n-1$ is equivalent to the game on $P_{n-5}^1$.  It holds $\gglls(P_n^2) \leq \gglls(P_{n-5}^1) + 2$ which is at most $\theta(P_n^1)$ in all the cases of $n \text{ (mod 5)}$.
\item Staller chooses $v_1=n-2$ or $v_1=n-3$: then Dominator can choose $v_2=n-3$ in the first case, and $v_2=n-2$ in the second. The game on the obtained predominated graph with predominated vertices $0, 2, 4, n-4, n-3, n-2$ and $n-1$ is equivalent to the game on $P_{n-3}^2$.  It holds $\gglls(P_n^2) \leq \gglls(P_{n-3}^2) + 2$, which is at most $\theta(P_n^2)$ in all the cases of $n \text{ (mod 5)}$.
\item Staller chooses $v_1=n-1$: then Dominator can choose $v_2=2$. The game on the obtained predominated graph with predominated vertices $0, 1, 2,3,4, n-2$ and $n-1$ is equivalent to the game on $P_{n-5}^1$.  It holds $\gglls(P_n^2) \leq \gglls(P_{n-5}^1) + 2$, which is at most $\theta(P_n^2)$ in all the cases of $n \text{ (mod 5)}$.
\item Staller chooses $4 < v_1 < n-3$: Then if Dominator chooses either $v_2=v_1+1$ or $v_2=v_1-1$, the obtained predominated graph has predominated vertices either $0, 2, 4, v_1-1, v_1, v_1+1, v_1+2, n-1$, or $0, 2, 4, v_1-2, v_1-1, v_1, v_1+1, n-1$ (if $v_1$ is 5 then consider only the first case, and notice that in cases $v_1$  is 5, 6 or 7  some vertices are written twice). Dominator can consider the same strategy as playing on two disjoint graphs: $P_{v_1}^2$ and $P_{n-(v_1+2)}^1$ in the first case or $P_{v_1-1}^2$ or $P_{n-(v_1+1)^1}$ in the second. We have:
\begin{eqnarray*}
\gglls(P_n^2) & \leq & \min\{\theta (P_{v_1}^2) + \theta (P_{n-(v_1+2)}^1) + 2, \theta (P_{v_1-1}^2) + \theta (P_{n-(v_1+1)}^1) + 2\} \\
& \leq & \theta (P_n^2)\,.
\end{eqnarray*}
In fact, the last inequality holds since $\theta (P_n^2) < \theta (P_{v_1}^2) + \theta (P_{n-(v_1+2)}^1) + 2$ only if $v_1 = 4 \text{ (mod 5)}$ and $n-(v_1+2) = 0 \text{ (mod 5)}$ as it can be checked by an easy examination. In that case the second entry of the minimum is smaller. 
\end{itemize}
This proves the assertion for $P_n^1$ and $P_n^2$.

Now we prove that $\ggll(P_n)\leq  \frac{4n}{5} + c_n^2$ by defining a strategy for Dominator. Let $v_1,\ldots,v_i$ be a sequence of moves. If $v_1,\ldots,v_i$ is also a legal sequence on $P_n^1$, then Dominator selects the same vertex as he would if the game was played on $P_n^1$. If the game is finished after $i$ moves on $P_n^1$, then Dominator can play an arbitrary move. If $v_{i}$ is an illegal move on $P_n^1$, then Dominator can play an arbitrary move. If $v_1,\ldots,v_i$ has some illegal moves (besides the last one) for the game on $P_n^1$, say $v_j$, then Dominator proceeds as if the sequence $v_1,\ldots,v_i$ is in fact $v_1,\ldots, v_{j-1}, v_{j+2},\ldots,v_i$. Since only two vertices in $P_n^1$ are predominated, the above procedure ensures that there are at most four moves more needed than on $P_n^1$.
\qed

\section{Problems, conjectures, and related extremal examples}
\label{sec:concluding}

We first demonstrate that the hierarchy of Theorem~\ref{thm:hierarchy} collapses for some graphs. For this sake recall that the Cartesian product $G\cp H$ of graphs $G$ and $H$ has the vertex set $V(G)\times V(H)$, vertices $(g,h)$ and $(g',h')$ being adjacent if either $gg'\in E(G)$ and $h=h'$, or $g=g'$ and $hh'\in E(H)$. 

\begin{prp}
\label{prop:Cartesian}
If $G$ is a connected graph with $n(G)\ge 2$ and $k\ge 2n(G)$, then 
$$\ggz(G\cp K_{1,k}) = \ggll(G\cp K_{1,k}) = 2n(G) - 1\,.$$ 
\end{prp}

\proof
Set $H = G\cp K_{1,k}$. Note that $\gamma(H) = \gamma_t(H) = n(G)$. Combining Theorem~\ref{thm:hierarchy} and Proposition~\ref{prp:bounds}(iii) we get
$$\ggz(H) \le \ggll(H)\le 2\gamma_t(H) - 1 = 2n(G) - 1\,.$$
Hence it remains to prove that $\ggz(H) \ge  2n(G) - 1$. The strategy of Staller is the following. Note that after the $i^{\rm th}$ move of Dominator, $i < n(G)$, he dominates at most $i$ subgraphs of $H$ induced by the set $V_u = \{(u,h):\ h\in V(K_{1,k})\}$, where $u\in V(G)$. These subgraphs are isomorphic to $K_{1,k}$ and are called fibers. Staller in each move follows Dominator in one of the fibers induced by $V_u$ in which Dominator has already played and for which there exists a neighbor $v$ of $u$ in $G$ such that Dominator did not yet play vertices from $V_v$. Clearly, Staller can select a vertex from $V_u$ which has a neighbor in $V_v$ that has not been dominated in the previous  moves. (Here we use the fact that $k\ge 2n(G)$.) Therefore the move of Staller is legal in the Z-domination game and hence the Z-domination game will last at least $2n(G) - 1$ moves. 
\qed

In view of Proposition~\ref{prop:Cartesian} we pose: 

\begin{prob} 
Characterize graphs $G$ without isolated vertices for which $\ggz(G) = \ggll(G)$ holds.
\end{prob} 

\noindent
Moreover we also pose: 

\begin{conj} 
\label{conj:trees}
If $T$ is a tree with $n(T)\ge 2$, then $\ggz(G) < \ggll(G)$ holds. 
\end{conj} 

\noindent
We have verified by computer that Conjecture~\ref{conj:trees} holds true for all trees on up to and including $18$ vertices. 

Recall from the end of Section~\ref{sec:hierarchy} that $\ggz(C_5) = \ggl(C_5)$ which implies that the other two sandwiched game domination parameters are also equal to $\ggz(C_5)$. This leads to: 

\begin{prob} 
Characterize graphs $G$ without isolated vertices for which $\ggz(G) = \ggl(G)$ holds.
\end{prob} 

Related to the examples presented in Section~\ref{sec:hierarchy} we also pose:

\begin{prob}
\label{prob:2n+1}
Is it true that $\ggll(G) \le 2\ggz(G) +1$ holds for an arbitrary graph $G$ without isolated vertices?  
\end{prob}

Note that from Proposition~\ref{prp:bounds} we easily get that $\ggll(G) \le 4\ggz(G) - 1$. Note also that if the answer to Problem~\ref{prob:2n+1} is affirmative, then the bound is best possible as demonstrated by any graph that contains a universal vertex. 

With respect to Section~\ref{sec:bounds} it would be interesting to systematically consider sharpness of the proved bounds and to characterize the graphs attaining the bounds. 

If $G$ is an isolate-free graph such that not all of its components are $K_2$, then by Theorem~\ref{thm:llbound} we have $\ggll(G)\le n(G)$. So it would be interesting to characterize the graphs that attain the equality. Instead, we propose the following special case which still seems very demanding. 

\begin{prob}
Characterize the trees $T$ with $\ggll(T) = n(T)$. 
\end{prob}

The $3/5$-conjecture for the domination game~\cite{KWZ-2013} and the $3/4$-conjecture for the total domination game~\cite{HKR-2017} are among the main sources of interest for the games, cf.~\cite{BKKR-2013, Bu-2015, Bu-2015b, HK-2016, HL-2017, Sc-2016} and~\cite{bujtas-2018, bujtas-2016, henning-2016}, respectively. An analogous question can be posed for the L-domination game hence we pose: 

\begin{conj}
If $G$ is a graph without isolated vertices, then $\ggl(G) \le \frac{6}{7}n(G)$. 
\end{conj}

The conjecture has been verified by computer for all trees up to $18$ vertices. It turned out that among them there are only two trees that attain the equality. These two trees belong to an infinite family of graphs $G$ for which $\ggl(G) = \frac{6}{7}n(G)$ holds which is defined as follows. Let $Y = S(K_{1,3})$, that is, $Y$ is the graph obtained from $K_{1,3}$ by subdividing each of its edges exactly once. If $G$ is an arbitrary graph, then let $G^Y$ be the graph obtained from $G$ by identifying each vertex of $G$ with the central vertex of a private copy of $Y$. In particular, the two trees mentioned above that were found by computer are $K_1^Y$ and $K_2^Y$.     

\begin{prp}
If $G$ is a graph, then $\ggl(G^Y) = \frac{6}{7}n(G^Y)$.  
\end{prp}

\proof
Let $G$ be an arbitrary graph with the vertex set $V(G)=\{u_1, \ldots, u_k\}$. To obtain $G^Y$, we attach a copy $Y^i$ of $Y$ to every $u_i$. The vertices of $Y^i$ are denoted by $v_1^i, \dots, v_7^i$ such that $v_1^i=u_i$ is the central vertex, $v_2^i$, $v_3^i$, $v_4^i$ are the support vertices, and $v_5^i$, $v_6^i$, $v_7^i$ are the leaves in $Y^i$. Clearly, $|V(G^Y)|=7k$.

First, we prove that $\ggl(G^Y) \le 6k$ and $\ggls(G^Y) \le 6k$. Consider the following strategy of Dominator.
\begin{itemize}
\item If it is a D-game, Dominator plays his first move in an arbitrary $Y^i$. 
\item Whenever Staller plays a vertex in a $Y^i$, Dominator replies with a move in the same $Y^i$, if there is such a legal move. Otherwise, Dominator may play in any $Y^j$ where a legal move can be made. 
\item Inside any $Y^i$, Dominator first plays the central vertex $v_1^i$ (if it has not been played by Staller earlier). Dominator's second move is a support vertex whose leaf neighbor has not been played yet. The third move may be any vertex.
\end{itemize}
By the first two rules, Dominator plays at least two of the first four moves in each $Y^i$. Hence, he can achieve that the central vertex and a support vertex are played before the adjacent leaf would be selected. These ensure that at least one leaf of $Y^i$ will not be played in the game. As a consequence, $\ggl(G^Y) \le 6k$ and $\ggls(G^Y) \le 6k$.

To prove the other direction, we first note that every support vertex of $G^Y$ must be played in the L-domination game. Hence, if Staller ensures that either all the three leaves or two leaves and the central vertex are played from every $Y^i$, then at most one vertex remains unplayed from each copy and therefore, the length of the game is at least $6n$. Consider the following strategy of Staller:
\begin{itemize}
\item If it is an S-game, Staller plays her first move in an arbitrary $Y^i$. 
\item Whenever Dominator plays a vertex in a $Y^i$, Staller replies in the same $Y^i$, if it is possible. Otherwise, she may choose any $Y^j$ where the game is not finished yet.
\item Inside any $Y^i$, Staller plays leaves while it is possible. Otherwise, she can choose any legal move from $Y^i$.
\end{itemize}
The first two rules ensure that Staller plays at least two of the first four moves in each $Y^i$. Since playing an (unplayed) leaf is not legal only if its support vertex and also the central vertex have been played earlier, the third rule ensures that at least three of the four vertices from $v_1^i$, $v_5^i$, $v_6^i$ and $v_7^i$ are played in the L-domination game. Since, as already noted, each of the $v_2^i$, $v_3^i$, and $v_4^i$ must be played, this proves that $\ggl(G^Y) \ge 6k$ and $\ggls(G^Y) \ge 6k$ and finishes the proof of the proposition.    
\qed


\end{document}